\documentclass{amsart}

\usepackage{graphicx}
\usepackage{xcolor}
\usepackage{amssymb}
\usepackage{mathrsfs}
\usepackage{xcolor}
\usepackage{tikz}
\usepackage{braids}
\usepackage{hyperref}
\usepackage{cite}

\newtheorem{theorem}{Theorem}[section]
\newtheorem{lemma}[theorem]{Lemma}
\newtheorem{corollary}[theorem]{Corollary}
\newtheorem{proposition}[theorem]{Proposition}

\newtheorem{conjecture}[theorem]{Conjecture}

\theoremstyle{definition}
\newtheorem{definition}[theorem]{Definition}
\newtheorem{example}[theorem]{Example}

\theoremstyle{remark}
\newtheorem{remark}[theorem]{Remark}

\numberwithin{equation}{section}

\begin{document}

\title{Translation distance bounds for fibered 3-manifolds with boundary}

\author{Alexander Stas}
\address{Mathematics Program, The Graduate Center, City University of New York, New York, NY}
\email{astas@gradcenter.cuny.edu}

\begin{abstract}
Given $M_\varphi$, a fibered 3-manifold with boundary, we show that the translation distance of the monodromy $\varphi$ can be bounded above by the complexity of an essential surface with non-zero slope. Furthermore we prove that the minimal complexity of a surface with non-zero slope in $M_{\varphi^n}$ tends to infinity as $n\to\infty$. Additionally, we show that an infinite family of fibered hyperbolic knots has translation distance bounded above by two, satisfying a conjecture by Schleimer which postulates that this behavior should hold for all fibered knots.
\end{abstract}

\maketitle

\section{Introduction}

Using geometry, combinatorics, and dynamics to study the topology of hyperbolic 3-manifolds is an active area of research. A theorem relating dynamics on a surface and fibered hyperbolic 3-manifolds is the following well known result of Thurston:

\begin{theorem}\cite{Thurston2-1998}
\label{paGivesFibered}
Let $\varphi:F\to F$ be a diffeomorphsim of a compact, connected surface with associated mapping torus $M_\varphi$. Then $M_\varphi$ is hyperbolic if and only if $\varphi$ is pseudo-Anosov.
\end{theorem}

Since periodic diffeomorphisms give rise to mapping tori that are Seifert manifolds and reducible diffeomorphisms give rise to mapping tori that have a non-trivial JSJ decomposition into simpler pieces, this result of Thurston implies that, in some sense, the only interesting fibered 3-manifolds are hyperbolic. 

A natural question to ask is ``what does $\varphi$ tell me about properties of $M_\varphi$?'' and vice-a-versa. The object of focus in this paper is the {\it translation distance} of $\varphi$ in the arc complex, denoted $d_\mathcal{A}(\varphi)$ (defined in Section \ref{defs}). Intuitively $d_\mathcal{A}(\varphi)$ is the minimum distance any vertex is moved in the arc complex under the action of $\varphi$.

In \cite{SY-2010}, Saito and Yamamoto show that any knot constructed by plumbing Hopf bands to the unknot, which were shown to be fibered by Harer in \cite{Harer-1982}, has translation distance at most 2 in the arc complex of the fiber, i.e. $d_\mathcal{A}(\varphi)\le2$. Unfortunately, their method doesn't determine how translation distance is affected by deplumbing. 

Bachman and Schleimer demonstrated a deep connection between incompressible surfaces in $M_\varphi$ and bounds on the translation distance in the curve complex, in particular they showed

\begin{theorem}\cite{BS-2005}
Let $\varphi:F\to F$ be a surface diffeomorphism of a closed orientable surface with genus $g(F)\ge2$. If $S\subset M_\varphi$ is a connected, orientable, incompressible surface, then either
\newpage
\begin{enumerate}
\item $S$ is isotopic to a fiber, or
\item $S$ is a torus and $d_\mathcal{C}(\varphi)\le1$, or
\item $d_\mathcal{C}(\varphi)\le|\chi(S)|$.
\end{enumerate}
\end{theorem}

The main objective of this paper is to extend this result of Bachman and Schleimer to fibered 3-manifolds with boundary. One nice application of this result will be to the class of fibered knot and link complements. To do so, we make use of the action of $\varphi$ on the arc and curve complex of the fiber. We prove

\begin{theorem}
\label{arcAndCurveTheorem}
Let $\varphi:F\to F$ be a diffeomorphism of a connected, compact, orientable surface $F$ with boundary and $\chi(F)<0$. Let $M_\varphi$ be the associated mapping torus. If $S$ is a connected, orientable, essential closed surface with no accidental parabolics or an essential non-longitudinal surface in $M_\varphi$, then either 
\begin{enumerate}
\item $\varphi$ is periodic and $d_\mathcal{AC}(\varphi^k)=0$ where $\varphi^k=\text{id}_F$,
\item  $\varphi$ is reducible and $d_\mathcal{AC}(\varphi)\le1$, or
\item $\varphi$ is pseudo-Anosov and $d_\mathcal{AC}(\varphi)\le|\chi(S)|$
\end{enumerate}
where $d_\mathcal{AC}(\varphi)$ is the translation distance of $\varphi$ in the arc and curve complex of the fiber.
\end{theorem}

If we restrict our attention further to the action of $\varphi$ on the arc complex, then the natural inclusion of the arc complex into the arc and curve complex gives a stronger result in the case where $\varphi$ is pseudo-Anosov. 

\begin{theorem}
\label{mainTheorem}
Let $\varphi:F\to F$ be a pseudo-Anosov diffeomorphism of a connected, compact, orientable surface $F$ with boundary and $\chi(F)<0$. Let $M_\varphi$ be the associated mapping torus. If $S$ is a connected, orientable, essential non-longitudinal surface properly embedded in $M_\varphi$, then $d_\mathcal{A}(\varphi)\le|\chi(S)|$ where $d_\mathcal{A}(\varphi)$ is the translation distance of $\varphi$ in the arc complex of the fiber.
\end{theorem}

Given a non-sporadic surface $F$ with boundary (i.e. $2g+b\ge 5$ where $g$ is the genus of $F$ and $b$ is the number of boundary components) and any $n\in\mathbb{Z}$, there are infinitely many homeomorphisms $\varphi:F\to F$ such that $d_\mathcal{A}(\varphi)>n$ (see Lemma \ref{MMforArc}). However, Schleimer has conjectured the following:

\begin{conjecture}[Schleimer \cite{Thompson-2016}]
\label{schleimerConj}
For any closed connected oriented 3-manifold $M$, there is a constant $t(M)$ with the following property: if $K \subset M$ is a fibered knot then the monodromy of $K$ has translation distance in the arc complex of the fiber at most $t(M)$. Furthermore $t(S^3) = 2$.
\end{conjecture}

Using Theorem \ref{mainTheorem} we give an infinite family of fibered hyperbolic knots in $S^3$ which satisfy Conjecture \ref{schleimerConj}. In particular we show 

\begin{corollary}
\label{infFamily}
Let $K\subset S^3$ be a fibered Montesinos knot with $r$--rational tangles and monodromy $\varphi$. If $r\ge 4$, then $d_{\mathcal{A}}(\varphi)\le 2$.
\end{corollary}

Define the {\it surface complexity of $M$} as $\Psi(M)=\min_{S\in\mathcal{S}}|\chi(S)|$ where $\mathcal{S}$ is the collection of all orientable non-longitudinal surfaces $S$ properly embedded in $M$. In Section \ref{surfaceComplexity}, we give an application of Theorem \ref{mainTheorem} and a result of Masur and Minsky \cite{MM1-1999} to show that essential surfaces become increasingly complex in $M_{\varphi^n}$ as $n$ increases. More specifically we show 

\begin{corollary}
\label{unboundedComplexity}
Let $F$ be a non-sporadic surface with boundary and $\varphi:F\to F$ a pseudo-Anosov diffeomorphism. Then the surface complexity $\Psi(M_{\varphi^n})\to\infty$ as $n\to\infty$.
\end{corollary}

\subsection{Acknowledgments}
We thank David Futer for mentioning this interesting problem and Saul Schleimer for his helpful conversations. We also thank Abhijit Champanerkar and Ilya Kofman for their countless suggestions and revisions. Lastly thanks to Daniel Berlyne, Alice Kwon, and Jacob Russell for many helpful discussions.


\section{Definitions and Background}
\label{defs}

For the entirety of this paper, we assume all 3-manifolds and surfaces in question are orientable. Let $F$ be a connected, compact, orientable surface with boundary and $\chi(F)<0$. Let $\varphi:F\to F$ be a diffeomorphism and $M_\varphi=F\times[0,1]/\sim$ the resultant mapping torus where $(x,0)\sim(\varphi(x),1)$. We say that a surface $S$ properly embedded in $M_\varphi$ is {\it essential} if $S$ is incompressible, boundary-incompressible, and not boundary parallel. If $\partial S=\emptyset$, than $S$ is essential if $S$ is incompressible and not boundary parallel. For example, each fiber $F(t) = F\times\{t\}$ is an essential surface. 

For what follows, it will be necessary to differentiate between surfaces that meet $\partial M_\varphi$ in curves parallel to $\partial F(t)$ and those that do not. Let $T^2_i$ be a boundary torus of $M_\varphi$. We define the longitude $\ell_i$ of $T^2_i$ to be the isotopy class of any component of $\partial F$ meeting $T^2_i$. Choose the meridian $m_i$ to be any simple closed curve on $T^2_i$ such that $\pi_1(T^2_i)=\langle m_i,\ell_i\rangle$. 

\begin{definition}
Let $T^2_i$ be a boundary torus of $M_\varphi$ and let $m_i,\ell_i$ be the generators of $\pi_1(T^2_i)$ defined above. Furthermore, let $S$ be an essential surface properly embedded in $M_\varphi$ which meets $T^2_i$. We define the {\it boundary slope of $S$ on $T^2_i$} to be the ratio $\frac pq\in\mathbb{Q}\cup\{\frac10\}$ where $[\partial S]=p[m_i]+q[\ell_i]\in H_1(T^2_i;\mathbb{Z})$. If $S$ does not meet $T^2_i$, the boundary slope of $S$ on $T^2_i$ is undefined. In this case we say that $S$ has slope $u$ on $T^2_i$. If $M_\varphi$ has $k$--boundary components, we can think of the slope of $S$ as the $k$--tuple $\overline{\frac pq}\in(\mathbb{Q}\cup\{\frac10\}\cup\{u\})^k$.
\end{definition}

The main theorem utilizes the boundary slope of an essential surface $S$ in $M_\varphi$ in several places. One of the key necessities for the proof is that in each boundary component of $M_\varphi$, the surface $S$ does not have slope $\frac01$. Thus, rather than stating this condition repeatedly, we make a definition.

\begin{definition}
Let $M_\varphi$ have $k$--boundary components. We say that a properly embedded surface $S\subset M_\varphi$ with slope $\overline{\frac pq}\in(\mathbb{Q}\cup\{\frac10\}-\frac01)^k$ is called a {\it non-longitudinal surface}. Note that a non-longitudinal surface meets every boundary component of $M_\varphi$.
\end{definition}

Although it is not known in general if every fibered 3-manifold contains such a surface, it was shown by Culler and Shalen in \cite{CS-1984} that all compact, connected, orientable 3-manifolds $M$ with $\partial M=T^2$ contain a properly embedded essential non-longitudinal surface proving a conjecture of Neuwirth.

In order to define the translation distance of the monodromy of the fibration, we need a metric space that $\varphi$ acts on naturally. Since $F$ has boundary, the complexes we will consider are the arc complex and the arc and curve complex of the fiber.

\begin{definition}
Given an essential arc (or curve) $\alpha$ on $F$ let $[\alpha]$ denote its isotopy class. Any collection of distinct isotopy classes of arcs (or arcs and curves) $A=\{[\alpha_0],\ldots,[\alpha_k]\}$ determines a $k$--simplex if for all $0\le i,j\le k$ with $i\neq j$ there are representatives $\alpha_i'\in[\alpha_i]$ and $\alpha_j'\in[\alpha_j]$ such that $\alpha_i'\cap\alpha_j'=\emptyset$. The {\it arc complex} $\mathcal{A}(F)$ (or {\it arc and curve complex} $\mathcal{AC}(F)$) is the simplicial complex determined by the union of all such simplices.
\end{definition}

\begin{remark}
We will refer to the isotopy class $[\alpha]$ as $\alpha$ unless otherwise stated. Furthermore, we are only concerned with the 0-- and 1--skeleta $\mathcal{A}^{(0)}(F)\subset \mathcal{A}^{(1)}(F)$ of $\mathcal{A}(F)$. 
\end{remark}

We can turn $\mathcal{A}^{(1)}(F)$ into a metric space by assigning length 1 to each edge and defining $d_\mathcal{A}(\alpha,\beta)$ to be the minimum path length over all paths from $\alpha$ to $\beta$. With this metric $\mathcal{A}^{(1)}(F)$ is $\delta$--hyperbolic, just as the curve complex of Masur and Minsky is $\delta$--hyperbolic \cite{MS-2013}. We turn $\mathcal{AC}^{(1)}(F)$ into a metric space in the same way and note that $\mathcal{AC}^{(1)}(F)$ is also $\delta$--hyperbolic \cite{KP-2010}. Note that the natural inclusion $\iota:\mathcal{A}^{(1)}(F)\to \mathcal{AC}^{(1)}(F)$ is distance-decreasing, i.e. $d_\mathcal{A}(\alpha,\beta)\ge d_\mathcal{AC}(\alpha,\beta)$ for any arcs $\alpha,\beta\subset F$.

As mentioned in the introduction, the object of interest is the translation distance of $\varphi$. Since $\varphi$ is a diffeomorphism, it acts on $\mathcal{A}^{(1)}(F)$ (resp. $\mathcal{AC}^{(1)}(F)$) as an isometry.

\begin{definition}
Let $\varphi:F\to F$ be a surface diffeomorphism. The {\it translation distance of $\varphi$ in} $\mathcal{A}^{(1)}(F)$ is
\[d_\mathcal{A}(\varphi)=\min\{d_\mathcal{A}(\alpha,\varphi(\alpha))~|~\alpha\in\mathcal{A}^{(0)}(F)\}.\]
The {\it translation distance of $\varphi$ in} $\mathcal{AC}^{(1)}(F)$ is defined analogously.
\end{definition}

The goal is to obtain a collection of arcs arising from the intersection of the fibers with a non-longitudinal surface which will define a path between $\alpha$ and $\varphi(\alpha)$. If done carefully, the Euler characteristic of the non-longitudinal surface can be used to bound the number of arcs in such a collection.


\section{Main Theorem}

For everything that follows, $F$ is a connected, compact, orientable surface with boundary and $\chi(F)<0$, and $\varphi:F\to F$ is a surface diffeomorphism. We start by proving a few necessary lemmas which enable us to prove all arcs and curves of intersection between the fibers and the essential surface are essential on both surfaces.

\begin{figure}
  \includegraphics[width=0.5\linewidth]{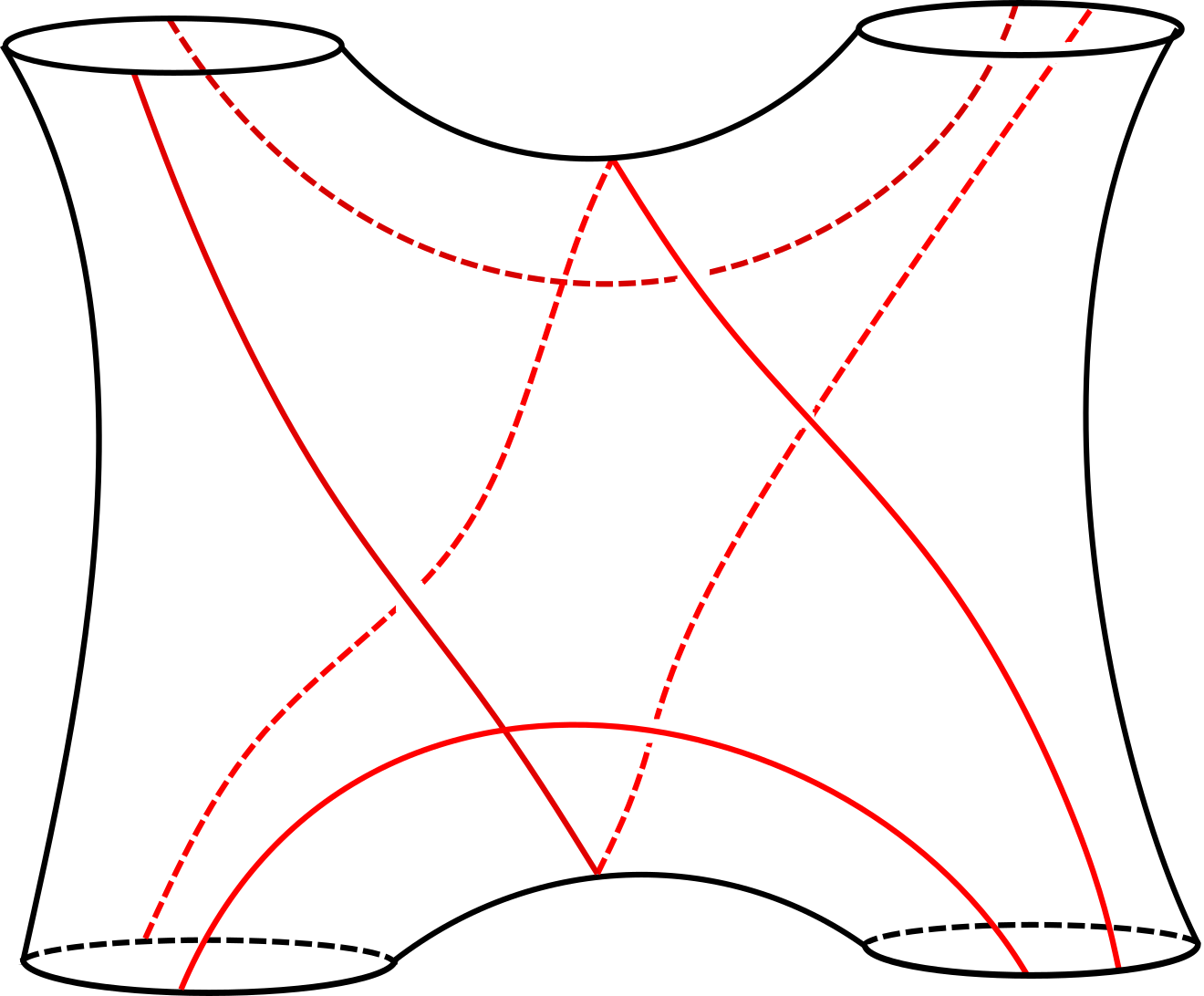}
  \caption{An induced singular foliation on the embedded Conway sphere, $C$, after isotopy. The red arcs are the non-transverse intersections of $C$ with the fibers. In particular there are exactly 2 intersections as $|\chi(C)|=2$.}
  \label{fig:conSphereFol}
\end{figure}

Suppose $M_\varphi$ is a fibered hyperbolic 3-manifold and let $\rho:\pi_1M_\varphi\to\text{PSL}(2,\mathbb{C})$ be a discrete faithful representation of the fundamental group of $M_\varphi$. Given $S$, a properly embedded surface in $M_\varphi$ or closed embedded surface, and an essential curve $\gamma\subset S$, we will say that $\gamma$ is an {\it accidental parabolic for $S$} if $\rho([\gamma])$ is a parabolic element under the induced representation. Therefore $\gamma$ is freely homotopic into $\partial M_\varphi$, but not boundary parallel in $S$. 

\begin{lemma}
Let $S$ be a connected, orientable, essential non-longitudinal surface properly embedded in a hyperbolic 3-manifold $M$. Suppose $\gamma$ is an accidental parabolic for $S$, let $\gamma'$ be the image of $\gamma$ in $T^2_i\subset\partial M$ under homotopy, and let $\partial S_i$ be the boundary components of $S$ in $T^2_i$. Then $\gamma'$ is parallel to $\partial S$ in $T^2_i$, i.e. $\iota(\gamma',\partial S_i)=0$.
\end{lemma}

\noindent{\bf Proof.} Let $A$ be the annulus arising from the homotopy of $\gamma$ into $T^2_i$, i.e. one boundary component of $A$ lies on $S$ and the other, $\gamma'$, in $\partial M$. Lemma 12.19 of Purcell \cite{Purcell-Book} produces an immersion of an annulus $A'$ into $M\backslash\backslash S$ such that one boundary component of $A'$ lies on a boundary component of a regular neighborhood of $S$ (possibly different from $\gamma$) and the other as before in $\partial M$. Since $A'$ is immersed in $M\backslash\backslash S$, we see that the boundary component of $A'$ in $\partial M$ must be parallel to $\partial S_i$. \qed \\

In \cite{ThurstonNorm}, Thurston showed that a surface $S$ in a fibered, compact 3-manifold with $\partial S$ contained in a fiber or transverse to the fibers can be (1) isotoped into a fiber, or (2) isotoped so that $S$ meets each fiber transversely except at finitely many singular points. Because we will use it later in a critical part of the proof, we prove a lemma quantifying the number of singularities arising from Thurston's construction.

\begin{lemma}
Suppose $S$ is a connected, orientable, essential closed surface with no accidental parabolics or an essential non-longitudinal surface properly embedded in $M_\varphi$. Then there exists an isotopy of $S$ such that
\begin{enumerate}
\item $S$ meets each fiber transversely except for at finitely many singular points,
\item each singular point contributes a 4-pronged singularity for the induced singular foliation of $S$, and
\item there are exactly $|\chi(S)|$ many such singular points, each contained in a unique fiber.
\end{enumerate}
\end{lemma}

\noindent{\bf Proof.} As in \cite{ThurstonNorm}, isotope $S$ so that it meets each fiber transversely except for at a finite number of saddle points. This proves (1) and (2). If necessary, we may perturb the isotopy slightly so that each singularity is contained in a unique fiber. By the Euler-Poincar\'e-Hopf formula, we have
\[2\chi(S)=\Sigma(2-P_s)\]
where the sum is taken over all singularities, $s$, in the induced singular foliation of $S$ and $P_s$ denotes the number of prongs at $s$. Since the induced foliation only contributes singularities with 4 prongs, there must be $|\chi(S)|$ many 4-pronged singularities (in the case $\partial S\neq\emptyset$, there are no boundary singularities as $\partial S$ meets the boundary of each fiber transversely). \qed\\

\begin{figure}
\begin{tikzpicture}
    \node[anchor=south west,inner sep=0] (image) at (0,0) {\includegraphics[width=0.7\textwidth]{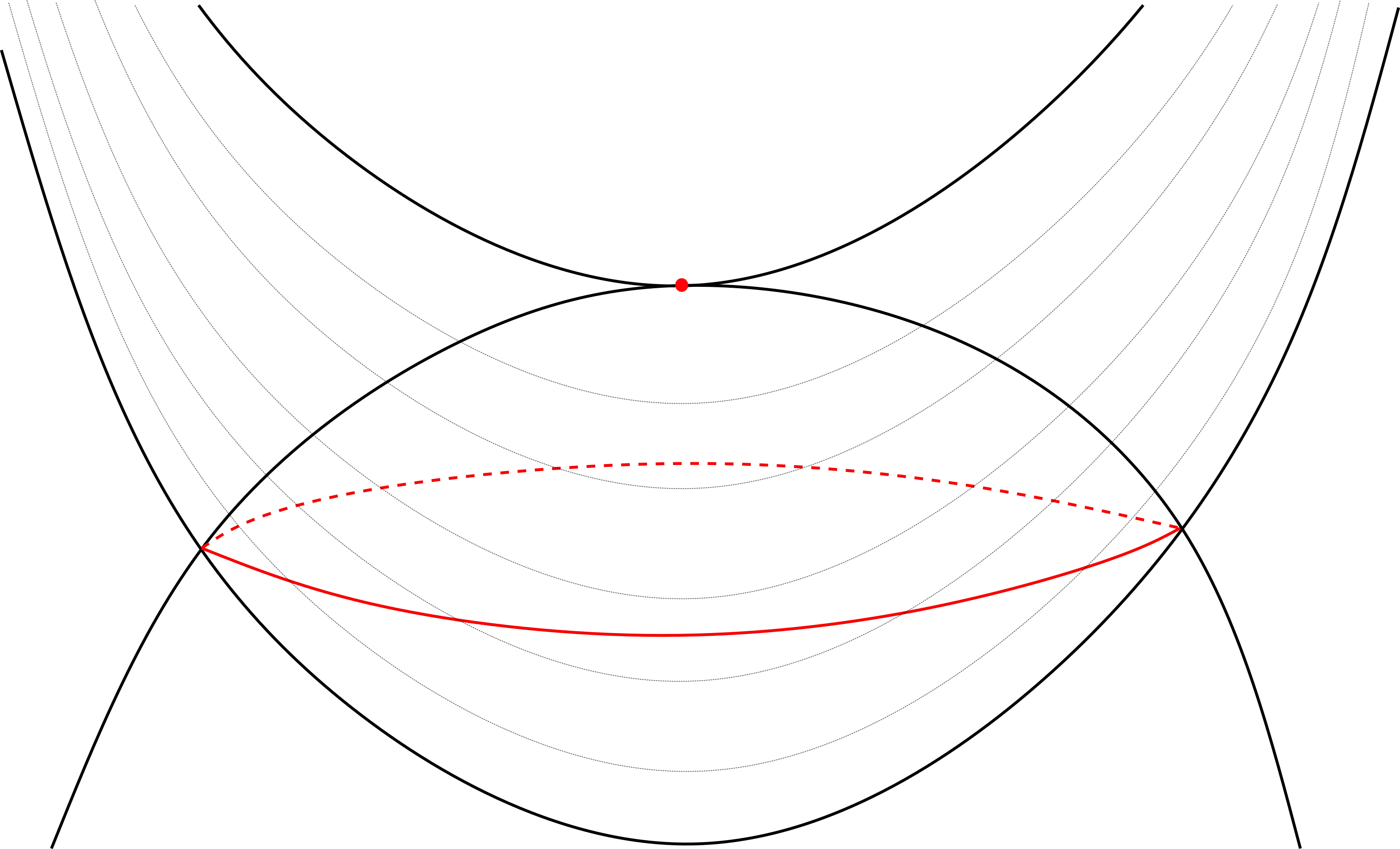}};
    \begin{scope}[x={(image.south east)},y={(image.north west)}]
        \draw (0.81,1.05) node {$F(t+\epsilon)$};
        \draw (1.0,1.05) node {$F(t)$};
        \draw (0.95,0.15) node {$S$};
        \draw (0.5,0.18) node {\textcolor{red}{$\gamma$}};
    \end{scope}
\end{tikzpicture}
      \caption{The red curve $\gamma$ bounds a disk on $S$ and $F(t)$. The resulting embedded 3--ball provides a center tangency (depicted in red) between $S$ and a fiber $F(t+\epsilon)$.}
    \label{fig:centerTangency}
\end{figure}
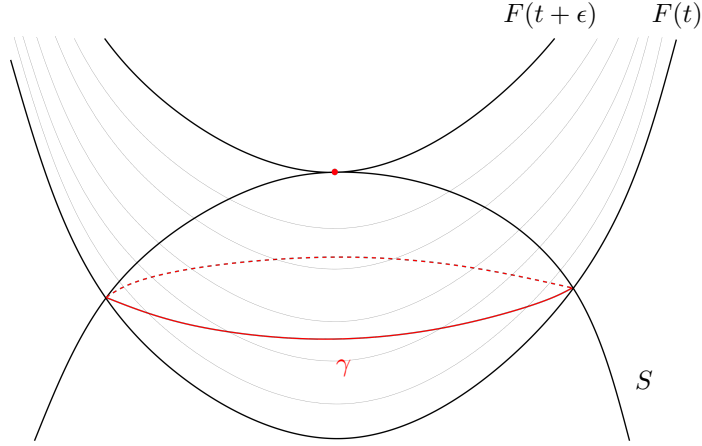

See Figure \ref{fig:conSphereFol} for a depiction of an induced singular foliation for a Conway sphere in a fibered knot complement.

\begin{lemma}
A curve $\gamma$ in a hyperbolic 3-manifold $M$ with boundary a collection of tori is homotopic into only one boundary component of $M$. 
\end{lemma}

\noindent{\bf Proof.} Suppose $\gamma$ can be homotoped into two distinct boundary components of $M$. The combination of these two homotopies gives an immersed annulus in $M$. The annulus theorem of Jaco \cite{Jaco-1980} provides an embedded annulus in $M$, contradicting the hyperbolicity of $M$.\qed\\

In the next two lemmas, we show that $S$ can be isotoped so that all curves and arcs of intersection with the fibers are essential on $S$ and the fibers.

\begin{lemma}
If $S$ is a connected, orientable, essential closed surface with no accidental parabolics or an essential non-longitudinal surface properly embedded in a fibered hyperbolic 3-manifold $M_\varphi$, then $S$ can be isotoped so that every curve of intersection between $S$ and any fiber is essential on both surfaces.
\end{lemma}

\noindent{\bf Proof.} Apply Lemma 3.2 to obtain an isotopy of $S$ such that $S$ is transverse to the fibers except at finitely many points and note that $S$ cannot be a fiber as $S$ is either non-longitudinal or closed. 

Let $F(t)=F\times\{t\}$ for $t\in[0,1]$ and consider $S\cap F(t)$ which is a collection of curves, arcs, or both. Our goal is to show that {\it every} transverse curve is essential on both surfaces. 

Suppose first that $\gamma$ is such a curve and $\gamma$ bounds a disk on $S$. Then $\gamma$ must also bound a disk on $F(t)$, for if it did not we would have found a compressing disk for $F(t)$, contradicting the incompressibility of the fibers. So $\gamma$ bounds a disk on both surfaces giving rise to an embedded 2--sphere. Since $M_\varphi$ is irreducible, this 2--sphere bounds a 3--ball $B$ which leads to a center tangency of $S$ and some fiber $F(t+\epsilon)$ contradicting transversality except at saddle points (see Figure \ref{fig:centerTangency}). A symmetric argument gives the case where $\gamma$ bounds a disk on $F(t)$.

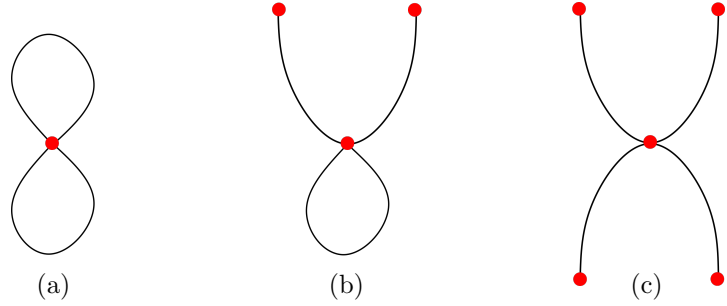
\begin{figure}
\begin{tikzpicture}
    \node[anchor=south west,inner sep=0] (image) at (0,0) {\includegraphics[width=0.75\textwidth]{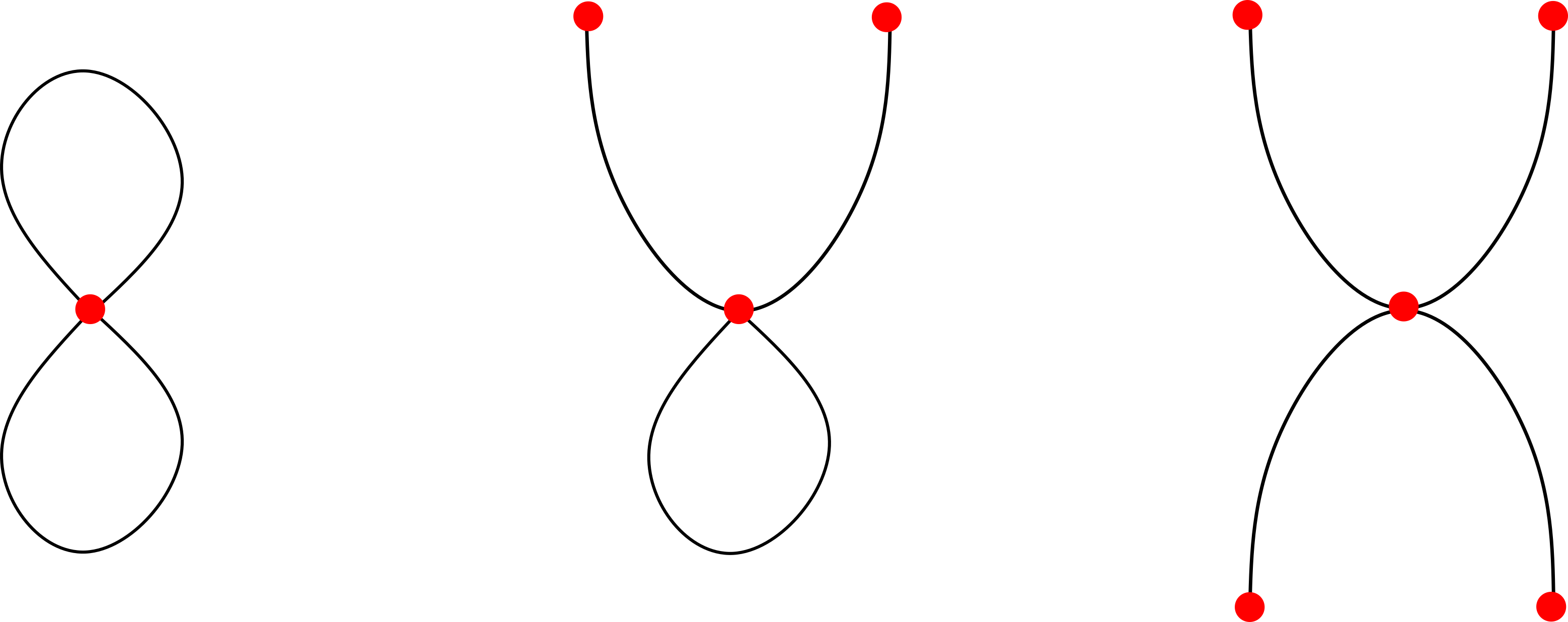}};
    \begin{scope}[x={(image.south east)},y={(image.north west)}]
        \draw (0.06,0.0) node {(a)};
        \draw (0.47,0.0) node {(b)};
        \draw (0.89,0.0) node {(c)};
    \end{scope}
\end{tikzpicture}
\caption{The possible types for singular components of $S\cap F(c_i)$ where (a) is the intersection of two curves, (b) is the intersection of a curve and an arc, and (c) is the intersection of two arcs.}
\label{fig:singComp}
\end{figure}

Suppose that $\partial S\neq\emptyset$ and $\gamma$ is isotopic into $\partial S$ and also isotopic into a component of $\partial F(t)$ which lies in a boundary torus $T_i^2$. Lemma 3.3 implies that such a curve $\gamma$ can only by homotopic into a single boundary component of $M_\varphi$ and so to avoid excessive notation, we simply call this boundary component $T^2$. The slope of $\partial F(t)$ in $T^2$ is $\frac01$ since $\partial F(t)$ is a longitude. By assumption, since $S$ is a non-longitudinal surface $\partial S\cap T^2$ is a collection of parallel non-longitudinal curves and thus $\partial S\cap T^2$ has non-zero slope. This leads to a contradiction since the combination of the above isotopies gives an isotopy between the longitude and a non-longitude of $T^2$ which is impossible.

Lastly, suppose that $\gamma$ is essential on $S$ and homotopic into $\partial F(t)$. If $S$ is a closed surface without accidental parabolics, then by definition such a curve $\gamma$ cannot exist as $\gamma$ is an accidental parabolic and the conclusion follows. Suppose therefore that $S$ is non-longitudinal. By Lemma 3.1 and Lemma 3.3 we have that $\gamma$ is homotopic to a curve $\gamma'$ in $T^2\subset\partial M_\varphi$ that is parallel to $\partial S\cap T^2$ and thus has the same non-zero slope as the non-longitudinal surface $S$. But by assumption $\gamma$ is homotopic to the longitude of $T^2$ as $\partial F(t)$ is a longitude. This gives a homotopy between a longitude and a non-longitude which is impossible. A symmetric argument gives the case where $\gamma$ is essential on $F(t)$ and homotopic into $\partial S$. Thus every transverse curve of intersection is essential on both surfaces.\qed

\begin{lemma}
If $S$ is a connected, orientable, non-longitudinal essential surface properly embedded in a fibered hyperbolic 3-manifold $M_\varphi$, then $S$ can be isotoped so that every arc of intersection between $S$ and any fiber is essential on both surfaces.
\end{lemma}

\noindent{\bf Proof.} Again, apply Lemma 3.2 to obtain an isotopy of $S$ such that $S$ is transverse to the fibers. In particular, we may assume we have the same isotopy as the one chosen in Lemma 3.4.

We now deal with the case where $F(t)$ and $S$ meet in an arc $\alpha$ inessential on $F(t)$ and essential on $S$. Then $\alpha$ being inessential on $F(t)$ provides us with a disk $D$ such that $\alpha$ cobounds $D$ with an arc $\beta\subset\partial F(t)\subset\partial M_\varphi$. Without loss of generality, we may assume that $\alpha$ is the innermost arc of intersection since if $\alpha'\neq\alpha$ is the innermost arc of intersection, $\alpha'$ is also inessential as $\alpha'\subset D$ and so $\alpha'$ also cobounds a disk.  By definition, $D$ is a boundary compression disk for $S$ contradicting the boundary-incompressibility of $S$. So $\alpha$ must be inessential on both surfaces. A symmetric argument shows that if $\alpha$ is inessential on $S$ then it must be inessential on $F(t)$.

Thus $\alpha$ must be inessential on $S$ and $F(t)$. In this case, $\alpha$ cobounds a disk $D$ on $S$ with $D\cap \partial S=\beta$, and $\partial D=\alpha\cup\beta$. Note that $\alpha$ is a leaf of the singular foliation of $S$ and $\beta$ must be met transversely at all points by our choice of isotopy of $S$. Also, the disk $D$ inherits a foliation which (possibly) contains finitely many finite-pronged singularities. Choose an arc $\alpha_1\subset D$ such that $\alpha_1$ cobounds a disk $D_1$ with a subarc $\beta_1\subset\beta$ and such that $\text{int}D_1$ does not contain any finite-pronged singularities. By identifying $D_1$ with a Euclidean disk, we may assume that the length of $\beta_1$ is 1. 

We will inductively construct a sequence of arcs $\{\beta_i\}$ as follows: Consider the disk $D_i$ with boundary $\alpha_i$ and $\beta_i$. Choose the midpoint $m_i\in\beta_i$. The arc $\alpha_{i+1}$ based at $m_i$ must hit another point $k_{i+1}\in\beta_i$. Let $\beta_{i+1}$ be the segment between $m_i$ and $k_{i+1}$ and so $\beta_{i+1}$ and $\alpha_{i+1}$ bound a disk $D_{i+1}\subset D_i$. Since $\ell(\beta_{i+1})<\frac{\ell(\beta_i)}{2}$ we have that $\ell(\beta_i)\to0$ as $i\to\infty$ and so there is a unique point $P$ that is the limit of $\{m_i\}$ and $\{k_i\}$.

Consider the arc $\alpha_P$ based at $P$ which must have its other endpoint $Q$ on $\beta$. Choose a small semi-circular neighborhood $N$ of $P$ not containing $Q$. Since $P$ is a limit point of $\{m_i\}$ and $\{k_i\}$, we can choose $j>0$ large enough so that $m_i,k_i\in N$ for all $i\ge j$. Therefore $\alpha_P$ must intersect $\alpha_{i+1}$ for all $i\ge j$ contradicting our choice of $D_1$ which did not contain any finite pronged singularities. Thus every arc of intersection must be essential on both surfaces.\qed\\

We now have the necessary tools to prove the main theorem. The proof proceeds by a similar argument to that of Bachman and Schleimer in their proof of Theorem 3.1 in \cite{BS-2005}. After proving Theorem \ref{mainTheorem}, we will use it to prove Theorem \ref{arcAndCurveTheorem}.\\

\noindent{\bf Theorem \ref{mainTheorem}.} {\it Let $\varphi:F\to F$ be a pseudo-Anosov diffeomorphism of a connected, compact, orientable surface $F$ with boundary and $\chi(F)<0$. Let $M_\varphi$ be the associated mapping torus. If $S$ is a connected, orientable, essential non-longitudinal surface properly embedded in $M_\varphi$, then $d_\mathcal{A}(\varphi)\le|\chi(S)|$ where $d_\mathcal{A}(\varphi)$ is the translation distance of $\varphi$ in the arc complex of the fiber.}\\

\noindent{\bf Proof.} Since $\varphi$ is pseudo-Anosov, we have by Theorem \ref{paGivesFibered} that $M_\varphi$ is hyperbolic. Use Lemma 3.2 to isotope $S$ to be transverse to the fibers except for at $n=|\chi(S)|$ singular points. Let $\{c_i\}_{i=0}^{n-1}$ be the values in $[0,1]$ where $F(c_i)$ fails to be transverse to $S$. Pick points $\{\tau_i\}_{i=0}^{n-1}$ in $[0,1]$ such that $c_{i-1}<\tau_i<c_i$ with indices taken mod $n$. Without loss of generality, we may assume $\tau_0=0$ by applying a rotation. 

\begin{figure}
\begin{tikzpicture}
    \node[anchor=south west,inner sep=0] (image) at (0,0) {\includegraphics[width=0.8\textwidth]{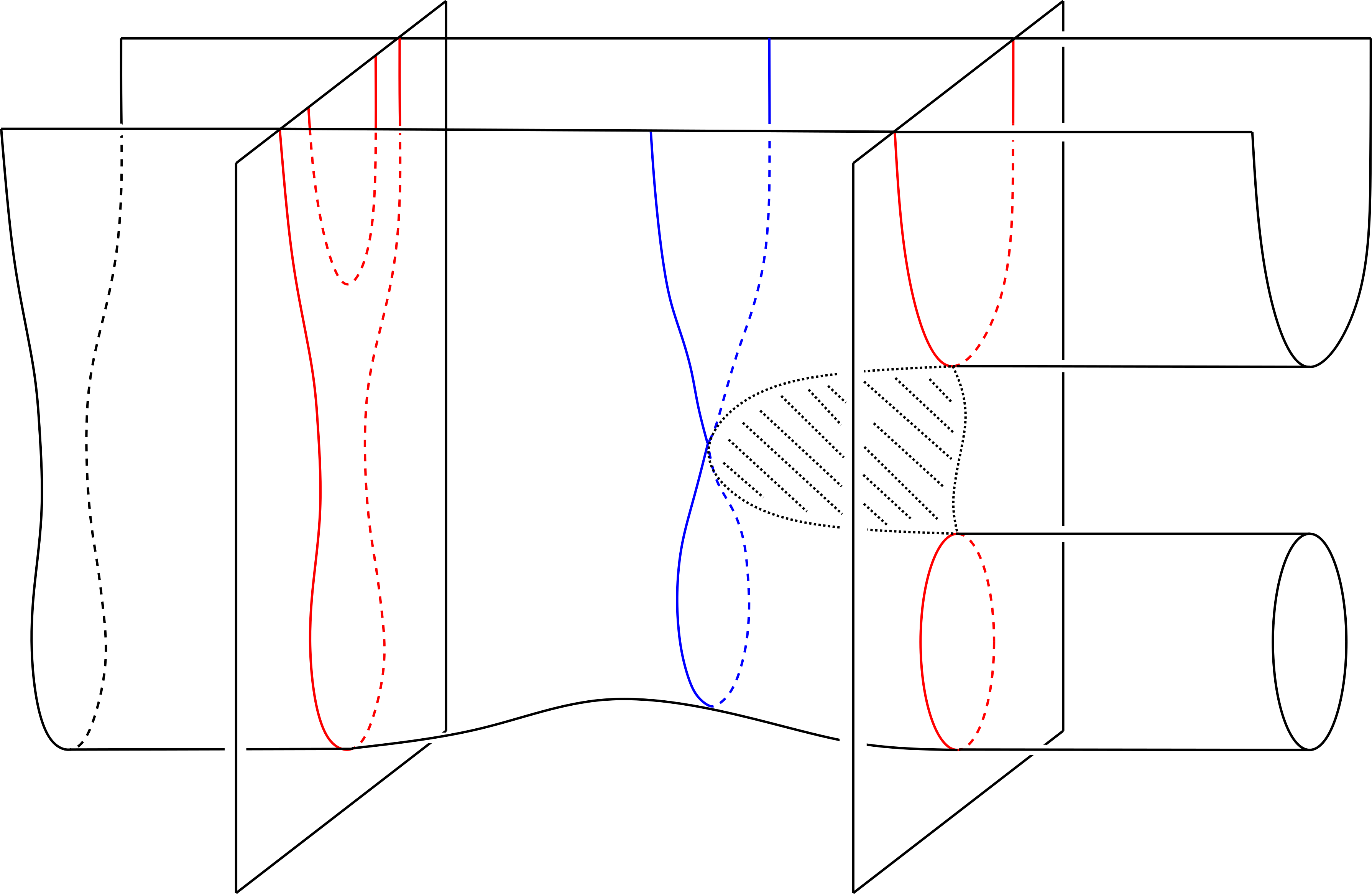}};
    \begin{scope}[x={(image.south east)},y={(image.north west)}]
        \draw (0.17,-0.05) node {$F(\tau_i)$};
        \draw (0.62,-0.05) node {$F(\tau_{i+1})$};
        \draw (0.73,1) node {$\alpha_{i+1}$};
        \draw (0.67,0.13) node {$\beta_{i+1}$};
        \draw (0.565,1.0) node {$\Gamma$};
        \draw (0.2,0.92) node {$\alpha_{i+1}''$};
        \draw (0.25,0.13) node {$\alpha_i$};
        \draw (0.73,0.5) node {$\gamma_1$};
        \draw (0.58,0.61) node {$\gamma_2$};
    \end{scope}
\end{tikzpicture}
\caption{A local picture of surface $S$ sitting in the product region $N=F\times[\tau_i,\tau_{i+1}]$. The result of isotoping $\alpha_{i+1}$ off of $F(\tau_{i+1})$ and into $F(\tau_i)$ is denoted $\alpha_{i+1}''$.}
\label{fig:annulusComp}
\end{figure}
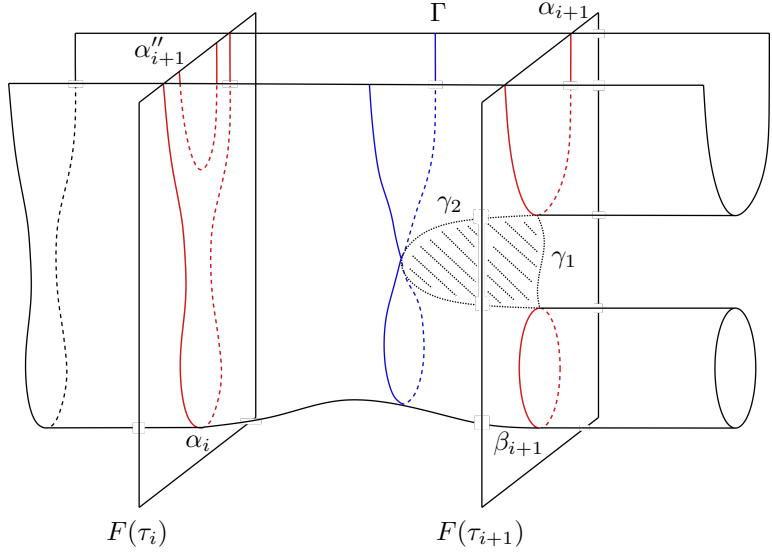

As the singular component of $S\cap F(c_i)$ must arise as one of three possibilities seen in Figure \ref{fig:singComp}, all singular points come from the intersection of two curves (Figure \ref{fig:singComp}(a)), an arc and a curve (Figure \ref{fig:singComp}(b)), or two arcs (Figure \ref{fig:singComp}(c)). 

We only need to consider case (b) and (c) since we will only be choosing an arc component for each $F(\tau_i)$. If $F(c_i)$ contains a singular component of type (a), then we may choose arcs $\alpha_{i-1}\subset S\cap F(\tau_{i-1})$ and $\alpha_i\subset S\cap F(\tau_i)$ such that $\alpha_{i-1}$ is isotopic to $\alpha_i$ through $S$. If this were not possible, then while pulling $\alpha_{i-1}$ forward into $F(\tau_i)$ through $S$ we would encounter a singular component of type (b) or (c) contradicting the assumption that the only singular component between $F(\tau_{i-1})$ and $F(\tau_i)$ is of type (a).

The proof in the remaining two cases is essentially the same and so we prove only one case. Since case (b) mixes arcs and curves, we prove (b).

Let $\Gamma$ be the singular component of $F(c_{i})\cap S$ which is a graph with three vertices and three edges. Let $N=F\times[\tau_{i},\tau_{i+1}]$ be the product region between our two non-critical values and let $P$ be the component of $S\cap N$ that contains $\Gamma$. We may assume that $P$ is the annulus in $N$ as seen in Figure \ref{fig:annulusComp} since by Morse theory, in a neighborhood of the singularity $P$ appears as a saddle. However, since our singular component arose from the intersection of an arc and a curve, $P$ is a saddle with a rectangle attached to the side of the saddle coming from the curve.

Therefore
\begin{enumerate}
\item $P$ is properly embedded in $N$,
\item $P\cap F(\tau_{i+1})$ is a disjoint essential arc $\alpha_{i+1}$ and essential curve $\beta_{i+1}$ in $F(\tau_{i+1})$, and
\item $P\cap F(\tau_i)$ is an essential arc $\alpha_i$ in $F(\tau_i)$.
\end{enumerate}

Now, $P$ is incompressible in $N$ as any compression disk for $P$ would give a compression disk for $S$. However, $P$ is boundary-compressible in $N$ via the disk $D$ bounded by an arc $\gamma_1$ in $F(\tau_{i+1})$ joining $\alpha_{i+1}$ to $\beta_{i+1}$ and an arc $\gamma_2$ on $P$ joining $\alpha_{i+1}$ to $\beta_{i+1}$. After this compression, we obtain $P'$ which is a properly embedded disk in $N$ giving an isotopy of $\alpha_{i}$ forward into $F(\tau_{i+1})$. Thus $\alpha_{i+1}$ is isotopic in $N$ to $\alpha_{i+1}''\subset F(\tau_i)$ which is disjoint from $\alpha_i$ as seen in Figure \ref{fig:annulusComp}.

From this construction we get a collection of $n+1$ arcs $\{\alpha_i'\}_{i=0}^n$ lying on $F(0)$ with the following properties: 

\begin{enumerate}
\item $\alpha_i'$ is isotopic to $\alpha_i$ through $F\times[0,\tau_i]$, $\alpha_0'=\alpha_0$, and $\alpha_n'$ is obtained by isotoping $\alpha_0$ off of $F(1)$ through $F\times[0,1]$,
\item $\alpha_i'\cap\alpha_{i+1}'=\emptyset$ for $i\in\{0,\ldots,n-1\}$, and
\item $\varphi(\alpha_n')$ is isotopic to $\alpha_0'$ since, by construction, we have that $\alpha_n'$ is isotopic to $\varphi^{-1}(\alpha_0')$.
\end{enumerate}

Thus
\[d_{\mathcal{A}}(\varphi)\le d_{\mathcal{A}}(\varphi(\alpha_n'),\alpha_n')=d_{\mathcal{A}}(\alpha_0',\alpha_n')\le\sum_{i=1}^{n}d_{\mathcal{A}}(\alpha_{i-1}',\alpha_i')\le n=|\chi(S)|.\]\qed

\vspace{0.3cm}
\noindent{\bf Proof of Theorem \ref{arcAndCurveTheorem}.} Suppose $\varphi$ is periodic and so $\varphi^k=\text{id}_F$. Therefore $\varphi^k$ fixes every arc and curve and so $d_\mathcal{AC}(\varphi^k)=0$.

If $\varphi$ is reducible, then $\varphi$ permutes some disjoint collection of isotopy classes of closed curves. If $\gamma$ is such a curve, then either $\gamma=\varphi(\gamma)$ or $\gamma$ and $\varphi(\gamma)$ are disjoint. In either case $d_\mathcal{AC}(\varphi)\le 1$.

Lastly, suppose that $\varphi$ is pseudo-Anosov. If $S\subset M_\varphi$ is a non-longitudinal surface, then since $d_\mathcal{AC}(\varphi)\le d_\mathcal{A}(\varphi)$ we are done by Theorem \ref{mainTheorem}. If $S$ is a closed surface without accidental parabolics, then there exists an isotopy of $S$ by Lemma 3.2 and Lemma 3.4 such that all curves of intersection between $S$ and any fiber are essential on all surfaces. Following the construction in the proof of Theorem \ref{mainTheorem}, we can choose a collection of $n+1$ curves $\{\gamma_i\}_{i=0}^n$ (rather than arcs) lying on the intersection of the fibers with $S$ satisfying properties (1)--(3) above. This is possible since after we isotopy $S$ to be transverse to the fibers away from finitely many singular points, all singular points will be of type (a) as in Figure \ref{fig:singComp}. Therefore
\[d_{\mathcal{AC}}(\varphi)\le d_{\mathcal{AC}}(\varphi(\gamma_n),\gamma_n)=d_{\mathcal{AC}}(\gamma_0,\gamma_n)\le\sum_{i=1}^{n}d_{\mathcal{AC}}(\gamma_{i-1},\gamma_i)\le n=|\chi(S)|.\]\qed


\section{Applications}

There are several interesting applications of Theorem \ref{mainTheorem} that we shall discuss here.

\subsection{Infinite families of low translation distance knots}

As mentioned in the introduction, Schleimer has conjectured\\

\noindent{\bf Conjecture \ref{schleimerConj}.} {\it For any closed connected oriented 3-manifold $M$, there is a constant $t(M)$ with the following property: if $K \subset M$ is a fibered knot then the monodromy of $K$ has translation distance in the arc complex of the fiber at most $t(M)$. Furthermore $t(S^3) = 2$.}\\

We start by showing that a consequence of Theorem \ref{mainTheorem} is that infinitely many knots satisfy Conjecture \ref{schleimerConj}.\\

\noindent{\bf Corollary \ref{infFamily}.} {\it Let $K\subset S^3$ be a fibered Montesinos knot with $r$--rational tangles and monodromy $\varphi$. If $r\ge 4$, then $d_{\mathcal{A}}(\varphi)\le 2$.}\\

\begin{figure}
\definecolor{linkcolor0}{rgb}{0.0, 0.0, 0.0}
\begin{tikzpicture}[line width=1, line cap=round, line join=round]
  \begin{scope}[color=linkcolor0]
    \draw (0.72, 3.45) -- (0.43, 3.74) -- (0.79, 4.10);
    \draw (0.79, 4.10) -- (1.15, 4.46) -- (1.87, 4.46) -- (2.24, 4.10);
    \draw (2.24, 4.10) -- (2.60, 3.74);
    \draw (2.60, 3.74) -- (2.96, 3.38);
    \draw (2.96, 3.38) -- (3.32, 3.02) -- (3.03, 2.73);
    \draw (2.89, 2.58) -- (2.60, 2.29) -- (2.24, 2.66);
    \draw (2.24, 2.66) -- (1.87, 3.02) -- (2.16, 3.31);
    \draw (2.31, 3.45) -- (2.52, 3.67);
    \draw (2.67, 3.81) -- (2.89, 4.03);
    \draw (3.03, 4.17) -- (3.32, 4.46) -- (4.04, 4.46) -- (4.33, 4.17);
    \draw (4.47, 4.03) -- (4.69, 3.81);
    \draw (4.84, 3.67) -- (5.05, 3.45);
    \draw (5.20, 3.31) -- (5.49, 3.02) -- (5.13, 2.66);
    \draw (5.13, 2.66) -- (4.76, 2.29);
    \draw (4.76, 2.29) -- (4.40, 1.93);
    \draw (4.40, 1.93) -- (4.04, 1.57) -- (4.33, 1.28);
    \draw (4.47, 1.14) -- (4.76, 0.85) -- (5.13, 1.21);
    \draw (5.13, 1.21) -- (5.49, 1.57) -- (5.20, 1.86);
    \draw (5.05, 2.01) -- (4.84, 2.22);
    \draw (4.69, 2.37) -- (4.47, 2.58);
    \draw (4.33, 2.73) -- (4.04, 3.02) -- (4.40, 3.38);
    \draw (4.40, 3.38) -- (4.76, 3.74);
    \draw (4.76, 3.74) -- (5.13, 4.10);
    \draw (5.13, 4.10) -- (5.49, 4.46) -- (6.21, 4.46) -- (6.57, 4.10);
    \draw (6.57, 4.10) -- (6.93, 3.74);
    \draw (6.93, 3.74) -- (7.29, 3.38);
    \draw (7.29, 3.38) -- (7.65, 3.02) -- (7.37, 2.73);
    \draw (7.22, 2.58) -- (6.93, 2.29) -- (6.57, 2.66);
    \draw (6.57, 2.66) -- (6.21, 3.02) -- (6.50, 3.31);
    \draw (6.64, 3.45) -- (6.86, 3.67);
    \draw (7.00, 3.81) -- (7.22, 4.03);
    \draw (7.37, 4.17) -- (7.65, 4.46) -- (8.38, 4.46) -- (8.74, 4.10);
    \draw (8.74, 4.10) -- (9.10, 3.74);
    \draw (9.10, 3.74) -- (9.46, 3.38);
    \draw (9.46, 3.38) -- (9.82, 3.02) -- (9.53, 2.73);
    \draw (9.39, 2.58) -- (9.10, 2.29) -- (8.74, 2.66);
    \draw (8.74, 2.66) -- (8.38, 3.02) -- (8.67, 3.31);
    \draw (8.81, 3.45) -- (9.03, 3.67);
    \draw (9.17, 3.81) -- (9.39, 4.03);
    \draw (9.53, 4.17) -- (9.82, 4.46) -- (9.82, 5.19) -- (0.43, 5.19) -- (0.43, 4.46) -- (0.72, 4.17);
    \draw (0.86, 4.03) -- (1.15, 3.74) -- (0.79, 3.38);
    \draw (0.79, 3.38) -- (0.43, 3.02) -- (0.43, 0.13) -- (9.82, 0.13) -- (9.82, 2.29) -- (9.46, 2.66);
    \draw (9.46, 2.66) -- (9.10, 3.02);
    \draw (9.10, 3.02) -- (8.74, 3.38);
    \draw (8.74, 3.38) -- (8.38, 3.74) -- (8.67, 4.03);
    \draw (8.81, 4.17) -- (9.10, 4.46) -- (9.46, 4.10);
    \draw (9.46, 4.10) -- (9.82, 3.74) -- (9.53, 3.45);
    \draw (9.39, 3.31) -- (9.17, 3.09);
    \draw (9.03, 2.95) -- (8.81, 2.73);
    \draw (8.67, 2.58) -- (8.38, 2.29) -- (7.65, 2.29) -- (7.29, 2.66);
    \draw (7.29, 2.66) -- (6.93, 3.02);
    \draw (6.93, 3.02) -- (6.57, 3.38);
    \draw (6.57, 3.38) -- (6.21, 3.74) -- (6.50, 4.03);
    \draw (6.64, 4.17) -- (6.93, 4.46) -- (7.29, 4.10);
    \draw (7.29, 4.10) -- (7.65, 3.74) -- (7.37, 3.45);
    \draw (7.22, 3.31) -- (7.00, 3.09);
    \draw (6.86, 2.95) -- (6.64, 2.73);
    \draw (6.50, 2.58) -- (6.21, 2.29) -- (5.49, 0.85) -- (5.20, 1.14);
    \draw (5.05, 1.28) -- (4.84, 1.50);
    \draw (4.69, 1.64) -- (4.47, 1.86);
    \draw (4.33, 2.01) -- (4.04, 2.29) -- (4.40, 2.66);
    \draw (4.40, 2.66) -- (4.76, 3.02);
    \draw (4.76, 3.02) -- (5.13, 3.38);
    \draw (5.13, 3.38) -- (5.49, 3.74) -- (5.20, 4.03);
    \draw (5.05, 4.17) -- (4.76, 4.46) -- (4.40, 4.10);
    \draw (4.40, 4.10) -- (4.04, 3.74) -- (4.33, 3.45);
    \draw (4.47, 3.31) -- (4.69, 3.09);
    \draw (4.84, 2.95) -- (5.05, 2.73);
    \draw (5.20, 2.58) -- (5.49, 2.29) -- (5.13, 1.93);
    \draw (5.13, 1.93) -- (4.76, 1.57);
    \draw (4.76, 1.57) -- (4.40, 1.21);
    \draw (4.40, 1.21) -- (4.04, 0.85) -- (3.32, 2.29) -- (2.96, 2.66);
    \draw (2.96, 2.66) -- (2.60, 3.02);
    \draw (2.60, 3.02) -- (2.24, 3.38);
    \draw (2.24, 3.38) -- (1.87, 3.74) -- (2.16, 4.03);
    \draw (2.31, 4.17) -- (2.60, 4.46) -- (2.96, 4.10);
    \draw (2.96, 4.10) -- (3.32, 3.74) -- (3.03, 3.45);
    \draw (2.89, 3.31) -- (2.67, 3.09);
    \draw (2.52, 2.95) -- (2.31, 2.73);
    \draw (2.16, 2.58) -- (1.87, 2.29) -- (1.15, 3.02) -- (0.86, 3.31);
  \end{scope}
\end{tikzpicture}
\caption{A projection of the Montesinos knot $K=M\left(\frac12,-\frac{2}{3},\frac{2}{5},-\frac{2}{3},-\frac{2}{3}~\Big|~0\right)$.}
  \label{fig:montesinosExample}
\end{figure}
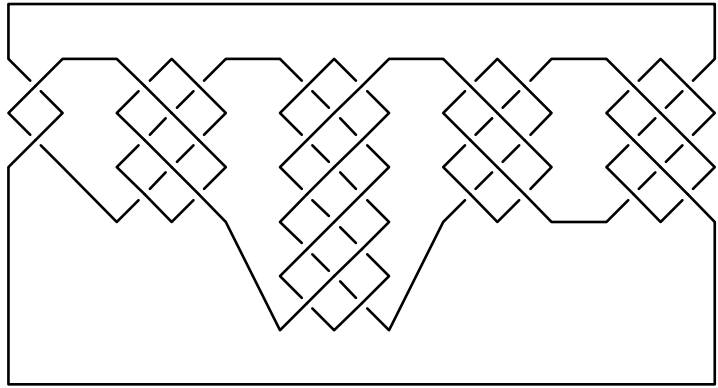

\noindent{\bf Proof.} In \cite{HM-2006}, Hirasawa and Murasugi determine when a Montesinos knot complement is fibered; in particular they show that there are infinitely many fibered Montesinos knots. Furthermore, Montesinos knots with $r\ge 4$ always contain an essential Conway sphere $S$ which has slope $\frac10$. Since $|\chi(S)|=2$, by Theorem \ref{mainTheorem} the conclusion follows.\qed

\begin{example}
Consider the Montesinos knot given by $K=M\left(\frac12,-\frac{2}{3},\frac{2}{5},-\frac{2}{3},-\frac{2}{3}~\Big|~0\right)$ using the notation of \cite{HM-2006}. Then $K$ is fibered hyperbolic, i.e. $S^3-K=M_\varphi$, and so by Corollary \ref{infFamily}, $d_{\mathcal{A}}(\varphi)\le 2$. See Figure \ref{fig:montesinosExample} for a projection of $K$.
\end{example}

Note the restriction in the above conjecture that $K$ is a knot and not a link. If we do not make that assumption, then the following construction (pointed out by Futer) provides infinitely many counterexamples:

\begin{proposition}
Conjecture \ref{schleimerConj} does not hold for links.
\end{proposition}
\noindent{\bf Proof.} Let $L$ be the braid closure of a psuedo-Anosov braid $\beta$ of $i$--components. Construct a link $L'$ of $(i+1)$--components from $L$ by adding the braid axis as a component. The link $L'$ is fibered with fiber surface the $i$--times punctured disk $D_i$ and monodromy given by the braid word $\beta$. We now get an infinite family of fibered links $\{L'_n\}_{n=1}^\infty$ where $L'_n$ is the braid closure of $\beta^n$ with the braid axis. By the above argument the fiber surface is $D_i$ and the monodromy is $\beta^n$. By Lemma \ref{MMforArc}, $d_\mathcal{A}(\beta^n)>kn-2$ for some constant $k$ depending only on $D_i$. Now choose $n$ large enough so that $kn-2>2$.\qed

\subsection{Lower bounds on surface complexity}
\label{surfaceComplexity}

The next application is of a more classical nature in the sense that we more closely rely on the action of the monodromy on the arc complex of the fiber. 

To prove the theorem of this subsection, we need to slightly extend a fundamental result of Masur and Minsky about translation distance in the curve complex which says

\begin{theorem}\cite{MM1-1999}
\label{masurMinsky}
For a non-sporadic surface $F$, there exists a constant $c>0$ such that for any pseudo-Anosov diffeomorphism $\varphi$ and any essential curve $\gamma$ on $F$, we have
\[d_\mathcal{C}(\gamma,\varphi^n(\gamma))\ge c|n|\]
for all $n\in\mathbb{Z}$.
\end{theorem}

We extend this result to the arc complex and obtain a similar lower bound for the distance between an arc and its iterates under some pseudo-Anosov map:

\begin{lemma}
\label{MMforArc}
For a non-sporadic surface $F$ with boundary, there exists a constant $k>0$ such that for any pseudo-Anosov diffeomorphism $\varphi$ and any essential {\it arc} $\alpha$ on $F$, we have
\[d_\mathcal{A}(\alpha,\varphi^n(\alpha))\ge k|n|-2\]
for all $n\in\mathbb{Z}$.
\end{lemma}

\noindent{\bf Proof.} In \cite{KP-2010}, Korkmaz and Papadopoulos showed that for any essential arc $\alpha$ on $F$, there is an essential curve $\gamma$ such that $\alpha$ and $\gamma$ are exactly distance 1 from each other in $\mathcal{AC}^{(1)}(F)$. Additionally, they showed that for any $x,y\in\mathcal{C}^{(0)}(F)$ we have that 
\[d_\mathcal{C}(x,y)\le 2d_\mathcal{AC}(x,y)\]

Let $\alpha$ be any essential arc on $F$ and choose $\gamma$ as above so that $d_\mathcal{AC}(\alpha,\gamma)=1$. By Theorem \ref{masurMinsky}, we have that
\begin{align*}
c|n|&\le d_\mathcal{C}(\gamma,\varphi^n(\gamma))\\
&\le 2d_\mathcal{AC}(\gamma,\varphi^n(\gamma))\\
&\le 2\left[d_\mathcal{AC}(\gamma,\alpha)+d_\mathcal{AC}(\alpha,\varphi^n(\alpha))+d_\mathcal{AC}(\varphi^n(\alpha),\varphi^n(\gamma))\right]\\
&=2\left[2+d_\mathcal{AC}(\alpha,\varphi^n(\alpha))\right].
\end{align*}
Letting $k=\frac c2$ we have
\[d_\mathcal{AC}(\alpha,\varphi^n(\alpha))\ge k|n|-2\]
and since $d_\mathcal{A}(\alpha,\varphi^n(\alpha))\ge d_\mathcal{AC}(\alpha,\varphi^n(\alpha))$ for any arc $\alpha$, the conclusion follows.\qed\\


Recall from the introduction the surface complexity of a hyperbolic 3-manifold $M$ is $\Psi(M)=\min_{S\in\mathcal{S}}|\chi(S)|$ where $\mathcal{S}$ is the collection of all non-longitudinal surfaces $S$ properly embedded in $M$. This is, in some sense, the ``right'' definition for the surface complexity, i.e. ruling out zero-slope surfaces and thus the fiber surface. If we include the Euler characteristic of the fiber surface, then $\Psi(M_{\varphi^n})=\Psi(M_{\varphi^m})$ for all $n,m\in\mathbb{Z}$ large enough. Additionally, since $M$ is hyperbolic, we have that $\Psi(M)>0$.\\

\noindent{\bf Proof of Corollary \ref{unboundedComplexity}.} For all $n\in\mathbb{Z}$, we have that $d_\mathcal{A}(\varphi^n)\ge k|n|-2$. As $n\to\infty$, the right-hand side gets arbitrarily large and thus the left-hand side does as well. By Theorem \ref{mainTheorem}, any essential non-longitudinal surface $S$ must satisfy $|\chi(S)|\ge k|n|-2$.\qed\\

\subsection{Primitiveness of knot monodromies and Schleimer's conjecture}\ \\
\hspace{0.5cm}We say that a mapping class $\varphi\in\text{Mod}(F)$ is {\it primitive} if whenever $\varphi=\psi^k$, we have that $\psi=\varphi$ and $k=1$. Intuitively this means that $\varphi$ is not a power of another mapping class. 

Let $\mathcal{K}_g$  be the collection of all fibered genus $g$ hyperbolic knots. It was proven by Stoimenow \cite{Stoimenow-2010} that there are infinitely many hyperbolic knots of genus $g\ge 2$. The question of whether or not there are infinitely many {\it fibered} hyperbolic knots of fixed genus is still open in general. However, it is believed that $|\mathcal{K}_g|=\infty$ for all $g\ge 2$ and some evidence in this direction is given by Kanenobu \cite{Kanenobu-1986} who shows that $|\mathcal{K}_2|=\infty$. The following result gives an interesting relationship between Schleimer's conjecture and the primitiveness of knot monodromies:

\begin{theorem}
Suppose that for every $n\in\mathbb{N}$ there exists $\psi\in\text{Mod}(F)$ which is primitive and a knot $K_n\in \mathcal{K}_2$ with monodromy given by $\varphi = \psi^m$ for some $m\ge n$, i.e. there exist knot monodromies of arbitrarily large powers. Then Conjecture \ref{schleimerConj} cannot hold.
\end{theorem}

\noindent{\bf Proof.} Choose $N$ large enough so that $kN-2>2$ where $k$ is as in Lemma \ref{MMforArc}. Then by assumption there exists a knot $K\in\mathcal{K}_2$ such that $S^3-K=M_{\psi^m}$ for some $m\ge N$ and $\psi$ primitive. Again, by Lemma \ref{MMforArc}, $d_\mathcal{A}(\psi^m)\ge km-2>2$ contradicting Conjecture \ref{schleimerConj}.\qed

\begin{remark}
Another way to rephrase this result is that if Schleimer's conjecture is to hold, there should be some uniform bound on how non-primitive a knot's monodromy can be, i.e. there exists $N\in\mathbb{N}$ such that if $\varphi=\psi^n$ is the monodromy of a knot and $\psi$ is primitive, then $n\le N$. 
\end{remark}

Mark Bell's program {\bf flipper} \cite{flipper} shows that the knots $7_7, 9_{48}, 12n0642, 12n0838$ are all fibered hyperbolic knots of genus $g=2$ with non-primitive monodromy.

\bibliographystyle{plain}
\bibliography{references}

\end{document}